\documentclass[12pt,reqno]{amsart}

\usepackage{amsfonts,color,amsthm,amsmath,amssymb}

\usepackage{color}
\def\Box{\vcenter{\vbox{\hrule\hbox{\vrule
     \vbox to 8.8pt{\hbox to 10pt{}\vfill}\vrule}\hrule}}}

\def\qed{{\hfill$\square$}}
\def\proof{{\vspace{-0.0cm}\bf Proof: \,}}

\def\Z{{\mathbb Z}}

\def\T{{\mathrm{Tr}}}

\def\F{{\mathbb F}}

\def\mod{{\mathrm{mod\,\,}}}

\def\Tr{{\mathrm{Tr}}}

\def\PG{{\mathrm{PG}}}
\def\AG{{\mathrm{AG}}}

\newtheorem{theorem}{Theorem}[section]

\newtheorem{lemma}[theorem]{Lemma}
\newtheorem{remark}[theorem]{Remark}

\newtheorem{proposition}[theorem]{Proposition}
\newtheorem{example}[theorem]{Example}

\numberwithin{equation}{section}

\begin{document}
\title[Conference matrices with maximum excess and  
two-intersection sets]
{Conference matrices with maximum excess and 
two-intersection sets}

\author[Koji Momihara and Sho Suda]{Koji Momihara$^*$ and  Sho Suda$^{\dagger}$}

\thanks{MSC2010: 11T22; 11T23; 05B20; 05E30}

\thanks{$^\ast$ Koji Momihara is supported by JSPS KAKENHI Grant Number (C)24540013.}
\thanks{$^{\dagger}$
Sho Suda is supported by  JSPS KAKENHI Grant Number 15K21075.}

\address{$^\ast$ Faculty of Education, Kumamoto University, 2-40-1 Kurokami, Kumamoto 860-8555, Japan} \email{momihara@educ.kumamoto-u.ac.jp}

\address{$^{\dagger}$ Department of Mathematics Education, Aichi University of 
Education, 1 Hirosawa, Igaya-cho, Kariya, Aichi 448-8542, Japan} \email{suda@auecc.aichi-edu.ac.jp}

\keywords{conference matrix; excess; quadratic residue; two-intersection set}

\begin{abstract}
A two-intersection set with parameters $(j;\alpha,\beta)$ for a block design is a $j$-subset of the point set of the  design, which intersects every block in $\alpha$ or $\beta$ points. 
In this paper, we show the existence of a two-intersection set with parameters $(2m^2-m+1;m^2-m,m^2)$ for the block 
design obtained from  translations of the set of nonzero squares in the finite field of order $q=4m^2+1$. 
As an application, we give a construction of conference matrices with maximum excess based on the two-intersection sets. 
\end{abstract}

\maketitle

\section{Introduction}\label{sec:1}
Let $P$ be a set of $v$ points and ${\mathcal B}$ be a collection of $b$ subsets of $P$, called {\it blocks}. We define $F=\{(p,B)\in P\times {\mathcal B}\,|\,p\in B\}$. Elements in $F$ are called {\it flags}.  The triple 
$(P,{\mathcal B},F)$ is called a {\it block design}. 
We say that $({\mathcal B},P,{\mathcal F}^\perp)$ with 
$F^\perp=\{(B,p)\,|\,(p,B)\in F\}$ is the {\it dual} of $(P,{\mathcal B},F)$. 
For convenience, we also say that the pair $(P,{\mathcal B})$ is  a block design. 

We consider a block design satisfying the following conditions: for each block $B\in {\mathcal B}$, there are exactly $k$ elements $p\in P$ such that $p\in B$. Dually, for each point $p\in P$, there are exactly $r$ blocks $B\in {\mathcal B}$ such that $p\in B$. Such a block design is called a {\it tactical configuration} or a {\it $1$-design}. It is clear that $vr=bk$. If for any two distinct 
points $a,b\in P$ the size of $\{B\in {\mathcal B}: a,b\in B\}$ is constant, say $\lambda$,  $(P,{\mathcal B})$ is called a {\it $2$-$(v,k,\lambda)$ design} or a {\it $2$-design}  for short. In particular, if $v=b$, it is called 
{\it symmetric}.  

Let $(P,{\mathcal B})$ be a block design and  $D$ be a $j$-subset of $P$. We say that $D$ is a 
{\it two-intersection set with parameters $(j;\alpha,\beta)$ for $(P,{\mathcal B})$} if the set  
\[
\{|B\cap D|\,:\,B\in {\mathcal B}\}
\]
contains exactly two numbers $\alpha$ and $\beta$. 
In this paper, we are interested in the existence of two-intersection sets and their applications. In particular, we consider a tactical configuration $(P,{\mathcal B})$  obtained from  translations of the set of nonzero squares in the finite field $\F_q$ of order $q\equiv 1\,(\mod{4})$, and in Section~\ref{sec:3}, 
we prove that there exists a two-intersection set with parameters $(2m^2-m+1;m^2-m,m^2)$ for $(P,{\mathcal B})$ if $q=4m^2+1$. 

Two-intersection sets have rich applications in algebraic combinatorics, in particular, for constructing strongly regular graphs and association schemes~\cite{BH,FMX,YBB} while there was no paper uniformly treating  two-intersection sets for block designs as far as the authors know.  In Section~\ref{sec:2}, we explain some of such applications briefly. Furthermore, we find a new application of two-intersection sets for constructing conference matrices  with maximum excess. 

A {\it conference matrix of order $n$} is an $n\times n$ $(0,-1,1)$-matrix $W$ with zero diagonal satisfying $WW^T=(n-1)I$, where $I$ is the $n\times n$ identity matrix. 
Conference matrices have been well-studied in relation to Hadamard matrices~\cite{IK}.  
Let $E(W)$ denote the sum of all entries of $W$. We say that $E(W)$ is the {\it excess} of $W$. We will show the following upper bound for excess of conference matrices in Appendix. 
\begin{proposition} \label{thm:confbound}
Let $W$ be a conference matrix of order $n$ and let $k$ be an odd integer such that  $k\le \sqrt{n-1}<k+2$. Then, it holds that
$E(W)\le \frac{n(k^2+2k+n-1)}{2(k+1)}$ with equality if and only if either one of the following holds: 
\begin{itemize}
\item[(i)] $n-1$ is a square and $W{\bf 1}_{n}=k{\bf 1}_n$; or 
\item[(ii)] $n-1$ is a nonsquare and $W{\bf 1}_{n}$ contains $k,k+2$ as its entries, 
\end{itemize}
where ${\bf 1}_n$ is the all one vector of length $n$. 
\end{proposition}
The excess of Hadamard matrices and complex Hadamard matrices have been studied in \cite{B,FK,K,KhS,KK,KS,S,XXS}. 
Note that regular conference matrices, that is conference matrices have the all-one vector as eigenvector, have the maximal excess, and see \cite{C} for constructions of regular conference matrices.
In this paper, we construct a conference matrix of order $n$ with $n-1$ nonsquare with excess attaining the upper bound of Proposition~\ref{thm:confbound} based on two-intersection sets obtained in Section~\ref{sec:3}. In particular, we will prove the following theorem. 
\begin{theorem}\label{thm:maxm}
For any prime power $q=p^r=4m^2+1$ with $p$ a prime congruent to $1$ modulo $4$, there exists a conference 
matrix of order $q+1$ with maximum excess attaining the bound of Proposition~\ref{thm:confbound}.   
\end{theorem}

\section{Two-intersection sets and their applications}\label{sec:2}
In this section, we consider  two-intersection sets for tactical configurations.  Let $(P,{\mathcal B})$ be a  
 tactical configuration and $D$ be a  two-intersection set with parameters $(j;\alpha,\beta)$ for $(P,{\mathcal B})$, i.e., $|D|=j$ and $\{|B\cap D|\,:\,B\in {\mathcal B}\}=\{\alpha,\beta\}$. It is clear  that the complement of 
$D$ is also a two-intersection set with parameters $(v-j;k-\alpha,k-\beta)$ for $(P,{\mathcal B})$. 
Set 
\[
D_\alpha^\perp=\{B\in {\mathcal B}\,:\,|B\cap D|=\alpha\}, \, \, D_\beta^\perp=
\{B\in {\mathcal B}\,:\,|B\cap D|=\beta\}. 
\]
Then, $|D_\alpha^\perp|+|D_\beta^\perp|=b$ and $\alpha |D_\alpha^\perp|+\beta |D_\beta^\perp|=
rj$. Hence, we have 
\begin{equation}\label{eq:bsize}
|D_\alpha^\perp|=\frac{-\beta b+rj}{\alpha-\beta}, \, \, 
|D_\beta^\perp|=b-\frac{-\beta b+rj}{\alpha-\beta}. 
\end{equation}
We say that each of $D_\alpha^\perp$ and $D_\beta^\perp$ is the {\it dual} of $D$. Here, the following  question naturally arises. Is the dual of $D$ also a two-intersection set for the dual of $(P,{\mathcal B})$? The answer is no in general, but there is a class of block designs giving an affirmative answer.     
\begin{proposition}
Let $(P,{\mathcal B})$ be a $2$-$(v,k,\lambda)$ design  and $D$ be a two-intersection set with parameters $(j;\alpha,\beta)$  
for $(P,{\mathcal B})$. Then, the dual $D_\alpha^\perp$ is also a two-intersection set 
with parameters $(j^\perp;\alpha^\perp,\beta^\perp)$ 
for the dual of $(P,{\mathcal B})$, where 
\[
j^\perp=\frac{rj-\beta b}{\alpha-\beta},\alpha^\perp=\frac{\lambda j-\beta r}{\alpha-\beta},  \beta^\perp=\frac{\lambda (j-1)+r-\beta r}{\alpha-\beta}. 
\] 
\end{proposition}
\proof 
Let $N$ be the matrix whose rows and columns are labeled by the elements of 
$P$ and ${\mathcal B}$, respectively, and entries are 
defined by 
\[
N_{p,B}=\begin{cases}
1& \text{ if } p\in B,\\
0& \text{ if } p \not\in B.  
\end{cases} 
\]
Let ${\bf x}$  and ${\bf y}$ be the $(0,1)$-vectors whose coordinates are  labeled by the elements of 
$P$ and ${\mathcal B}$, respectively, and entries are defined by 
\[
{\bf x}_{p}=\begin{cases}
1& \text{ if } p\in D,\\
0& \text{ if } p \not\in D,   
\end{cases} 
\mbox{\,  and \, }
{\bf y}_{B}=\begin{cases}
1& \text{ if } B\in D_\alpha^\perp,\\
0& \text{ if } B \in D_\beta^\perp.    
\end{cases} 
\]
It is clear that 
\begin{equation}\label{eq:both}
{\bf x}^TN=\alpha {\bf y}^T+\beta ({\bf 1}-{\bf y})^T. 
\end{equation}
By multiplying both sides of \eqref{eq:both} by $N^T$ from right, we have 
\[
{\bf x}^TNN^T=(\alpha-\beta) {\bf y}^TN^T+\beta {\bf 1}^TN^T. 
\]
Let $I$ be the identity matrix of order $v$ and $J$ be 
the all-one matrix of order $v$.  
Since $NN^T=\lambda J+(r-\lambda)I$ and $ {\bf 1}^TN^T=r{\bf 1}^T$, 
${\bf y}^TN^T$ has exactly two entries $\alpha^\perp=(\lambda j-\beta r)/(\alpha-\beta)$ and $\beta^\perp=(\lambda (j-1)+r-\beta r)/(\alpha-\beta)$.  Furthermore, 
by multiplying both sides of \eqref{eq:both} by ${\bf 1}$ from right, we have 
\[
{\bf x}^TN{\bf 1}=(\alpha-\beta) {\bf y}^T{\bf 1}+\beta {\bf 1}^T{\bf 1}. 
\]
Since ${\bf x}^TN{\bf 1}=r{\bf x}^T{\bf 1}=rj$ and $ \beta {\bf 1}^T{\bf 1}=\beta b$, we have $j^\perp={\bf y}^T{\bf 1}=(rj-\beta b)/(\alpha-\beta)$. 
Hence, the dual $D_\alpha^\perp$  is  a two-intersection set 
for the dual of $(P,{\mathcal B})$. 
\qed

\begin{remark}{\em 
If $(P,{\mathcal B})$ is a $2$-design, the parameter $j$ is determined by 
the other parameters. Let ${\mathcal B'}=\{B\cap D\,|\,B\in {\mathcal B}\}$. Double-counting the number of  pairs of distinct points of  $(D,{\mathcal B}')$, we have ${j\choose 2}\lambda={\alpha \choose 2}|D_\alpha^\perp|+{\beta \choose 2}|D_\beta^\perp|$. Substituting \eqref{eq:bsize}, the parameter $j$ is computable. In the language of design theory, $(D,{\mathcal B}')$ is a {\it pairwise balanced design} with two block sizes. 
}\end{remark}

Two-intersection sets have been studied in 
algebraic combinatorics in relation to strongly regular graphs and association schemes. Hereafter, we will assume that the reader is familiar with the theory of association schemes. 
\begin{example}\label{ex:apptwo1}
{\em ({\bf Projective two-intersection set}) \,
Let $\F_q$ be the finite field of order $q$ and $\F_{q}^\ast=\F_q\setminus \{0\}$. 
Let $(P,{\mathcal B})$ be the $2$-design  
obtained from
points and hyperplanes of the $n$-dimensional  projective space $\PG(n,q)$ over $\F_q$. 
Let $D$ be a two-intersection set for $(P,{\mathcal B})$. Define \[
C=\{x y: x \in D,y\in \F_q^\ast\}. 
\] 
We consider the graph $\Gamma=(V,E)$ defined as  $V=\F_{q}^{n+1}$ and $(x,y)\in E$ if and only  if $x-y \in C$, which is called a {\it Cayley graph} on 
$\F_{q}^{n+1}$. The set $C$ is called the {\it connection set} of  $\Gamma$. It is known that this Cayley graph $\Gamma$ forms a strongly regular graph~\cite[p.~134]{BH}. 
The set $D$ is particularly called a {\it projective two-intersection set}. The existence of projective  two-intersection sets has been well-studied 
in finite geometry. 
}   
\end{example}

\begin{example}\label{ex:apptwo2} {\em ({\bf Affine two-intersection set}) \,
Let $(P,{\mathcal B})$ be the $2$-design obtained from
points and hyperplanes of the $n$-dimensional  affine space $\AG(n,q)$ over $\F_q$. 
Let $D$ be a  two-intersection set for $(P,{\mathcal B})$. We now assume that $(P,{\mathcal B})$ is a {\it residual} of the $2$-design obtained from   points and hyperplanes of $\PG(n,q)$, i.e., $P$ is the complement of a fixed hyperplane $H$ of $\PG(n,q)$ and the blocks $B\in {\mathcal B}$ are the restrictions of hyperplanes of $\PG(n,q)$ to $P$. 
Then, we can regard $D$ as a subset of the set of projective points of 
$\PG(n,q)$. 
Then, the Cayley graphs on $\F_{q}^{n+1}$ with connection sets 
\[
C_1=\{x y\,:\, x \in D,y\in \F_q^\ast\}, \, C_2=\{x y\,:\, x \in P\setminus D,y\in \F_q^\ast\}, \, C_3=\{x y\,:\, x \in H,y\in \F_q^\ast\}
\]
partition the complete graph on $V=\F_{q}^{n+1}$. In particular, this partition forms a $3$-class association scheme~\cite{FMX}. The set $D$ is called an {\it affine two-intersection set}. A first infinite family of affine two-intersection sets was recently found in \cite{DDMR,FMX}. 
}  
\end{example}
\begin{example}\label{ex:apptwo3}{\em ({\bf Relative $2$-design})\, 
Let $(P,{\mathcal B})$ be a symmetric $2$-design.
Assume that there exists a two-intersection set $D$ for $(P,{\mathcal B})$ satisfying $|D|=k$ and $D\not \in {\mathcal B}$. 
Then, $(P,{\mathcal B})$ forms a tight relative $2$-design with respect to $D$ in the Johnson association scheme $J(v,k)$.
Some constructions and all possible parameters with $v\leq 100$ were given in \cite{YBB}. 
All the examples with $v\leq 100$ have the structure of coherent configurations. 
} 
\end{example}

We now give a new application of 
two-intersection sets for tactical configurations  obtained from quadratic residues of finite fields.  Let $q\equiv 1\,(\mod{4})$ be a prime power and 
$S$ be the set of nonzero squares of $\F_q$. Set $P=\F_q$ and 
\begin{equation}\label{eq:desquad}
{\mathcal B}=\{\{x+a:x \in S\}:a \in \F_q\}. 
\end{equation}
Then, $(P,{\mathcal B})$ is a tactical configuration with  $v=b=q$, $k=r=(q-1)/2$. Note that $(P,{\mathcal B})$ does not form a $2$-design. The set $S$ is called the 
{\it Paley partial difference set},  which satisfies that the list $\{x-y: x,y\in S,x\not=y\}$ covers every element of $S$ (resp. $\F_q^\ast\setminus S$) exactly $(q-5)/4$ times (resp.  $(q-1)/4$ times). It is well-known that the Cayley graph on $\F_q$ with connection set $S$ forms a strongly regular graph. The Paley partial difference sets also have an application for constructing conference matrices. 
Let $M$ be a $q\times q$ $(0,1,-1)$-matrix whose rows and columns are labeled by the elements of $\F_q$ and entries are defined by 
\[
M_{i,j}=\begin{cases}
0& \text{ if } j-i=0,\\
1& \text{ if } j-i\in S,\\
-1& \text{ if } j-i\in \F_q^\ast \setminus S. 
\end{cases} 
\]
Define 
\begin{equation}\label{eq:conf}
W=\begin{pmatrix} 0  & {\bf 1}_q^T \\ {\bf 1}_q & -M \end{pmatrix}. 
\end{equation}
Then, $W$ forms a conference matrix. We construct a  conference matrix with maximum excess  by switching the signs of some rows and columns of $W$. 
\begin{theorem}\label{thm:const}
Let $q=4m^2+1$ be a prime power and $(P,{\mathcal B})$ be the block design defined  in \eqref{eq:desquad}.  Assume that there is a two-intersection set with parameters $(2m^2-m+1;m^2-m,m^2)$ for $(P,{\mathcal B})$.  
Then, there exists a  conference 
matrix $W'$ of order $q+1$ such that $W'{\bf 1}_{q+1}$ has entries $2m-1$ and 
$2m+1$.  
\end{theorem}
\proof 
Let $D$ be the assumed two-intersection set and $W$ be the conference matrix defined  in \eqref{eq:conf}. Set $\alpha=m^2-m$ and $\beta=m^2$.  
Multiply by $-1$ the columns indexed by the elements of $D$ of 
$\begin{pmatrix} {\bf 1}_q^T \\ -M \end{pmatrix}$.   
Denote the resulting matrix by 
$\begin{pmatrix} {\bf b}^T \\ M' \end{pmatrix}$.  
Then, multiply by $-1$  the rows indexed by the elements of $D_{\alpha}^\perp$ of 
$\begin{pmatrix} {\bf 1}_q & M' \end{pmatrix}$.  Denote the resulting matrix by $\begin{pmatrix} {\bf c} & M'' \end{pmatrix}$.  Then, 
$W'=\begin{pmatrix} 0  & {\bf b}^T \\ {\bf c} & M'' \end{pmatrix}$ is the desired  conference matrix. 

It is clear that ${\bf b}^T{\bf 1}_{q}=2m-1$ by the assumption that $|D|=2m^2-m+1$. Furthermore, since $(\alpha,\beta)=(m^2-m,m^2)$, we have 
\[
\left(\begin{pmatrix} {\bf 1}_q & M' \end{pmatrix}{\bf 1}_{q+1}\right)_i=
\begin{cases}
2m-1& \text{ if } i\in D_{\beta}^\perp\mbox{ and  } i\not \in D,\\
2m+1& \text{ if } i\in D_{\beta}^\perp \mbox{ and  } i \in D,\\
-2m-1& \text{ if } i\in D_{\alpha}^\perp \mbox{ and  } i\not \in D,  \\
-2m+1& \text{ if } i\in D_{\alpha}^\perp \mbox{ and  } i\in D. 
\end{cases}
\]
Hence, $W'{\bf 1}_{q+1}$ has entries 
$2m-1,2m+1$. 
 \qed
\vspace{0.3cm}

In the next section, we will construct two-intersection sets satisfying the condition of Theorem~\ref{thm:const}. Then, Theorem~\ref{thm:maxm}  immediately follows by Proposition~\ref{thm:confbound}.   
\section{Construction of two-intersection sets}\label{sec:3}
\subsection{Preliminary on characters of finite fields}\label{sec:31}
In this section, we will assume that the reader is familiar with the basic theory of characters of finite fields. 

For a positive integer $m$, set $\zeta_m=\exp^{\frac{2\pi \sqrt{-1}}{m}}$. 
Let $q=p^r$ be a prime power with $p$ a prime. 
For a multiplicative character
$\chi$  and the canonical
additive character $\psi$ of $\F_q$, we define the {\it Gauss sum} by
\[
G_q(\chi)=\sum_{x\in \F_q^\ast}\chi(x)\psi(x) \in \Z[\zeta_{q-1},\zeta_p].
\]
We list a few basic properties of Gauss sums below:
\begin{enumerate}
\item[(i)] $G_q(\chi)\overline{G_q(\chi)}=q$ if $\chi$ is nontrivial;
\item[(ii)] $G_q(\chi^{-1})=\chi(-1)\overline{G_q(\chi)}$;
\item[(iii)] $G_q(\chi)=-1$ if $\chi$ is trivial.
\end{enumerate}
Let $\omega$ be a primitive element of $\F_q$ and $k$ be a positive integer dividing $q-1$. 
For $0\le i\le k-1$ we set $C_i^{(k,q)}=\omega^i C$, where $C$ is the multiplicative subgroup of index $k$ of $\F_q^\ast$. By the orthogonality of characters, the sums
$\psi(C_i^{(k,q)})=\sum_{x\in C_i^{(k,q)}}\psi(x)$, $0\le i\le k-1$, so-called  
{\it Gauss periods}, can be expressed as a linear combination of Gauss sums:
\begin{equation}
\psi(C_i^{(k,q)})=\frac{1}{k}\sum_{j=0}^{k-1}G_q(\chi^{j})\chi^{-j}(\gamma^i), \; 0\le i\le k-1,
\end{equation}
where $\chi$ is a multiplicative character of order $k$ of $\F_q$.  For example, if $k=2$, 
we have
\begin{equation}\label{eq:Gaussquad}
\psi(C_i^{(2,q)})=\frac{-1+(-1)^iG_q(\eta)}{2},\; 0\le i\le 1,
\end{equation}
where $\eta$ is the quadratic character of $\F_q$. In particular, the quadratic 
Gauss sum is explicitly computable. 
\begin{theorem}\cite[Theorem~5.15]{LN97} \label{thm:Gauss}
Let $q=p^s$ be a prime power with $p$ a prime and $\eta$ be the quadratic character of $\F_q$. 
Then, 
\begin{equation}\label{eq:Gaussquad1}
G_q(\eta)=\begin{cases}
(-1)^{s-1}q^{1/2}& \text{ if }  p\equiv 1\,(\mod{4}), \\
(-1)^{s-1}\zeta_4^s q^{1/2} & \text{ if } p\equiv 3\,(\mod{4}). 
\end{cases}
\end{equation}
\end{theorem} 
Furthermore, 
we need to define {\it Jacobi sums}. We extend the domain of multiplicative characters $\chi$ of $\F_q$ to all elements of $\F_q$ by setting 
$\chi(0)=1$ or $\chi(0)=0$ depending on whether $\chi$ is trivial or not. 
For multiplicative characters $\chi_1$  and $\chi_2$  of $\F_q$, define 
\[
J(\chi_1,\chi_2)=\sum_{x\in \F_q}\chi_1(x)\chi_2(1-x) \in \Z[\zeta_{q-1}].  
\]
In this paper, we treat Jacobi sums $J(\chi_1,\chi_2)$ with $\chi_1$ 
the quadratic character and $\chi_2$ a multiplicative  character of order $4$ of $\F_q$. 
\begin{lemma}\label{lem:facto}{\em (\cite{St})} 
Let $q\equiv 1\,(\mod{4})$ be  a prime power. 
Let $\eta$ be the quadratic character and $\chi$ a multiplicative  character of order $4$ of $\F_q$. Put $J(\eta,\chi)=a+b\zeta_4\in \Z[\zeta_4]$. Then, 
\begin{itemize}
\item[(i)] $a\equiv -1\,(\mod{4})$ if $q\equiv 1\,(\mod{8})$, 
\item[(ii)] $a\equiv 1\,(\mod{4})$ if $q\equiv 5\,(\mod{8})$.  
\end{itemize}
Conversely, for any prime power $q=p^r=a^2+b^2\equiv 1\,(\mod{4})$ with $p\equiv 1\,(\mod{4})$ a prime satisfying  (i) or (ii) above and $\gcd{(a,q)}=1$, it holds that $J(\eta,\chi)=a+b\zeta_4$, where the sign of $b$ is ambiguously determined. If $p\equiv 3\,(\mod{4})$, $r$ is even and  $J(\eta,\chi)=a$. 
\end{lemma}
In this paper, we do not need to care about the signs of $a,b$. 

We will use the following formula on Jacobi sums in the next section. 
\begin{proposition}\label{prop:Jaco}{\em (\cite[Exercise~5.60]{LN97})}
For any $a,b\in \F_q^\ast$ and a multiplicative character $\chi$ of $\F_q$, it holds that 
\[
\sum_{x\in \F_q}\chi(ax^n+b)=\chi(b)\sum_{j=1}^{d-1}{\chi'}^{-j}(a){\chi'}^j(-b)
J({\chi'}^j,\chi), 
\]
where $\chi'$ is a multiplicative character of order $d=\gcd{(n,q-1)}$ of $\F_q$. 
\end{proposition}
\subsection{Construction}\label{sec:32}
Let $q=p^r=4m^2+1$ be a prime power with $p$ a prime congruent to $1$ modulo $4$, and let 
$\omega$ be a primitive element of $\F_{q^2}$.  
Let $\chi_4$ be a multiplicative character of order $4$ of $\F_{q^2}$, and 
$\eta$ and $\chi_4'$ be multiplicative characters of order $2$ and 
$4$ of $\F_{q}$, respectively. Assume that 
$\chi_4(\omega)=\chi_4'(\omega^{q+1})=\zeta_4$. By Lemma~\ref{lem:facto}, there are $\epsilon,\delta\in \{-1,1\}$ such that $J(\eta,\chi_4')=\epsilon+2m\delta \zeta_4$.  

Let $\ell$ be an integer not divisible by $q+1$, and put $n=\omega^{\ell(q+1)}$ and $t=\omega^\ell+\omega^{\ell q}$. 
Fix $h\in \{0,1,2,3\}$ and $\ell$ so that the following conditions are satisfied: 
\begin{equation}\label{eq:cond2}
\chi_4'(n)=\zeta_4^{\frac{-\epsilon \delta+1}{2}+h}, \quad 
\chi_4'(n-t^2/4)=-\zeta_4^{\epsilon+2h}.  
\end{equation}
We will see in Remark~\ref{rem:cond2} that such a pair $(h,\ell)\in \{0,1,2,3\}\times \{0,1,\ldots,q^2-2\}$ always exists.

%
\begin{theorem}\label{thm:twoint}
Let $q=p^r=4m^2+1$ be a prime power with $p$ a prime congruent to $1$ modulo $4$, and let 
$\omega$ be a primitive element of $\F_{q^2}$.  Let $h$ and $\ell$ be integers 
defined as above. 
Define
\[
D_{\ell,h}=\left\{x\in \F_q\,|\, 1+x\omega^\ell\in C_h^{(4,q^2)} \cup C_{h+1}^{(4,q^2)} \right\}.  
\]
Then, the set $\{|D_{\ell,h}\cap (C_0^{(2,q)}+s)|:s\in \F_q\}$ contains exactly two numbers $m^2-m$ and $m^2$, and $|D_{\ell,h}|=2m^2-m+1$. 
\end{theorem}
This theorem implies that $D_{\ell,h}$ is a two-intersection with parameters $(j;\alpha,\beta)=(2m^2-m+1;m^2-m,m^2)$ for the block design $(P,{\mathcal B})$ defined in \eqref{eq:desquad}. 
We prove this theorem by a series of propositions below. 
\begin{proposition}\label{prop:charac1}
Let $q\equiv 1\,(\mod{4})$ be a prime power and $\chi_4$ be a multiplicative
character of order $4$ of $\F_{q^2}$. Let $\omega$ be a primitive element of $\F_{q^2}$.  
Put $n=\omega^{\ell(q+1)}$ and $t=\omega^\ell+\omega^{\ell q}$. Then, 
\[
\sum_{x\in \F_q}\chi_4(1+\omega^\ell x)={\chi_{4}'}^3(n)
{\chi_{4}'}^{3}(n-t^2/4)J(\eta,\chi_4'), 
\]
where $\eta$ and $\chi_4'$ are multiplicative characters of order $2$ and $4$ of $\F_q$ such that $\chi_4(\omega)=\chi_4'(\omega^{q+1})$.  
\end{proposition}
\proof
Let $\chi_8$ be a multiplicative character of order $8$ of $\F_{q^2}$ such that $\chi_8^{q+1}=\chi_4$. Note that the restriction of $\chi_{8}$ to $\F_q$ 
is of order $4$, which coincides with $\chi_4'$. In fact, 
\[
\chi_4'(\omega^{q+1})=\chi_4(\omega)=\chi_8^{q+1}(\omega)=\chi_8(\omega^{q+1}). 
\]
Then, we have 
\begin{align}
\sum_{x\in \F_q}\chi_{4}(1+\omega^\ell x)
=&\, \sum_{x\in \F_q}\chi_{8}((1+\omega^\ell x)^{q+1})\nonumber\\
=&\, \sum_{x\in \F_q}\chi_{4}'(1+t x+n x^2)\nonumber\\
=&\,\sum_{x\in \F_q}\chi_{4}'(nx^2+1-t^2/(4n)). 
\label{eq:com2}
\end{align}
By Proposition~\ref{prop:Jaco}, the summation \eqref{eq:com2} is reformulated as 
\begin{align*}
\sum_{x\in \F_q}\chi_{4}'(nx^2+1-t^2/(4n))=&\,\eta(n)
{\chi_{4}'}^{3}(1-t^2/(4n))J(\eta,\chi_4')\\
=&\,{\chi_{4}'}^3(n)
{\chi_{4}'}^{3}(n-t^2/4)J(\eta,\chi_4'). 
\end{align*}
This completes the proof. \qed

\begin{proposition}\label{prop:charac2}
With the notations of Proposition~\ref{prop:charac1}, assume that 
$\ell$ is not divisible by $q+1$. Then, 
for any $s\in \F_q$,
\[
\sum_{x\in \F_q\setminus \{s\}}\chi_{4}(1+\omega^\ell x)\eta(x-s)={\chi_4'}^3(u)
{\chi_{4}'}^{3}(n-t^2/4)J(\eta,\chi_4')-\chi_{4}'(n), 
\]
where $u=1+ts+ns^2$. 
\end{proposition}
\proof 
Since the restriction of $\chi_4$ to $\F_q$ is of 
order $2$, we have 
\begin{align}
\sum_{x\in \F_q\setminus \{s\}}\chi_{4}(1+\omega^\ell x)\eta(x-s)=&\, \sum_{y\in \F_q^\ast}\chi_{4}(1+\omega^\ell (y+s))\chi_4^{-1}(y)\nonumber\\
=&\, \sum_{y\in \F_q^\ast}\chi_{4}(y^{-1}(1+\omega^\ell s)+\omega^\ell)\nonumber\\
=&\, \sum_{y\in \F_q}\chi_{4}(y(1+\omega^\ell s)+\omega^\ell)-\chi_4(\omega^\ell). \label{eq:chara2}
\end{align}
Note that $\chi_4(\omega^\ell)=\chi_4'(n)$ and $1+\omega^\ell s\not=0$. 
Setting  $u=(1+\omega^\ell s)(1+\omega^{\ell q} s)=1+t s+n s^2(\not=0)$ and 
$v=w^{\ell q}(1+\omega^\ell s)+w^{\ell }(1+\omega^{\ell q} s)$, we have 
\begin{align}
\eqref{eq:chara2} =&\,\sum_{y\in \F_q}\chi_{8}((y(1+\omega^\ell s)+\omega^\ell)^{q+1})-\chi_4'(n)\nonumber\\
=&\, \sum_{y\in \F_q}\chi_{4}'(u y^2+v y+  n)-\chi_{4}'(n)\nonumber\\
=&\,\sum_{y\in \F_q}\chi_{4}'(uy^2+n-v^2/(4u))-\chi_{4}'(n). 
\label{eq:com3}
\end{align} 
By Proposition~\ref{prop:Jaco}, the summation of the left hand side of \eqref{eq:com3} is reformulated as 
\[
\sum_{x\in \F_q}\chi_{4}'(ux^2+n-v^2/(4u))=\eta(u){\chi_{4}'}^{3}(n-v^2/(4u))J(\eta,\chi_4'). 
\]
Noting that $n-v^2/(4u)=(4n-t^2)/(4u)$, we obtain 
\[
\eqref{eq:com3}={\chi_4'}^3(u)
{\chi_{4}'}^{3}(n-t^2/4)J(\eta,\chi_4')-\chi_{4}'(n). 
\]
This completes the proof. 
\qed
\vspace{0.3cm}

We are now ready for proving Theorem~\ref{thm:twoint}. 

{\bf  Proof of Theorem~\ref{thm:twoint}:}\,  
The characteristic functions of $C_0^{(2,q)}$ and 
 $C_h^{(4,q^2)}\cup C_{h+1}^{(4,q^2)}$ are, respectively,  given as 
\[
g(x)=\frac{1}{2}(\eta(x)+1), \quad  x\in \F_q^\ast,
\]
and 
\[
f(x)=\frac{1}{4}\sum_{j=h,h+1}\sum_{i=0}^3\zeta_4^{-ji}\chi_4^i(x), \quad  x\in \F_{q^2}^\ast. 
\]
The size $N_s$ of the set $D_{\ell,h}\cap (C_0^{(2,q)}+s)$ is expressed as 
\[
\sum_{x\in \F_q\setminus \{s\}}f(1+\omega^\ell x)g(x-s). 
\]
By the definitions of $g(x)$ and $f(x)$,  we have 
\begin{align}
N_s=&\,\frac{1}{8}\sum_{x\in \F_q\setminus \{s\}}(\eta(x-s)+1)\Big(\sum_{j=h,h+1}\sum_{i=0}^3\zeta_4^{-ji}\chi_4^i(1+\omega^\ell x)\Big)\nonumber\\
=&\,\frac{1}{8}\sum_{x\in \F_q\setminus \{s\}}(\eta(x-s)+1)\Big(2+\zeta_4^{-h}(1-\zeta_4)\chi_4(1+\omega^\ell x)+\zeta_4^{-3h}(1+\zeta_4)\chi_4^3(1+\omega^\ell x)\Big). \label{eq:com1}
\end{align}
Let $N_{s,1}=\sum_{x\in \F_q\setminus \{s\}}\chi_4(1+\omega^\ell x)$ and 
$N_{s,2}=\sum_{x\in \F_q\setminus \{s\}}\chi_4(1+\omega^\ell x)\eta(x-s)$. 
Then, 
\begin{equation}\label{eq:comN}
\eqref{eq:com1}=
\frac{1}{8}\Big(
2q-2+\zeta_4^{-h}(1-\zeta_4)(N_{s,1}+N_{s,2})+\overline{\zeta_4^{-h}(1-\zeta_4)(N_{s,1}+N_{s,2})}
\Big). 
\end{equation}
By Propositions~\ref{prop:charac1} and \ref{prop:charac2}, we have 
\begin{align}\label{eq:nnn}
N_{s,1}+N_{s,2}=({\chi_4'}^3(n)+{\chi_4'}^3(u)){\chi_4'}^3(n-t^2/4)J(\eta,\chi_4')-(\chi_4'(n)+\chi_4'(u)). 
\end{align}
Here, we used $\chi_4(1+\omega^\ell s)=\chi_4'(u)$. 
Substituting $\chi_4'(n)=\zeta_4^{\frac{-\epsilon \delta+1}{2}+h}$, $ 
\chi_4'(n-t^2/4)=-\zeta_4^{\epsilon+2h}$ and $J(\eta,\chi_4')=\epsilon+2m\delta \zeta_4$ into \eqref{eq:nnn} and continuing from \eqref{eq:comN}, we have
\begin{equation}\label{eq:dual}
N_s=
\begin{cases}
m^2-m& \text{ if $\epsilon \delta=-1$ and $\chi_4'(u)\in \{\zeta_4^{h+1},\zeta_4^{h+2}\}$} 
\\
& \text{\quad or  $\epsilon \delta=1$ and $\chi_4'(u)\in \{\zeta_4^{h},\zeta_4^{h+3}\}$,}\\
m^2& \text{ if  $\epsilon \delta=-1$ and $\chi_4'(u)\in \{\zeta_4^h,\zeta_4^{h+3}\}$ }\\
& \text{\quad or  $\epsilon \delta=1$ and $\chi_4'(u)\in \{\zeta_4^{h+1},\zeta_4^{h+2}\}$.}
\end{cases} 
\end{equation}
Thus,  $N_s$ takes exactly two values according to
$s$. 

Next,  we compute the size of $D_{\ell,h}$:   
\begin{equation}\label{eq:size}
|D_{\ell,h}|=\sum_{x\in \F_q}f(1+x \omega^\ell)=
\frac{1}{4}\sum_{x\in \F_q}\Big(2+\zeta_4^{-h}(1-\zeta_4)\chi_4(1+\omega^\ell x)+
\zeta_4^{-3h}(1+\zeta_4)\chi_4^3(1+\omega^\ell x)
\Big). 
\end{equation}
By Proposition~\ref{prop:charac1}, we have 
\begin{equation}\label{eq:sizech}
\sum_{x\in \F_q}\chi_4(1+\omega^\ell x)={\chi_{4}'}^3(n)
{\chi_{4}'}^{3}(n-t^2/4)J(\eta,\chi_4'). 
\end{equation}
Substituting $\chi_4'(n)=\zeta_4^{\frac{-\epsilon \delta+1}{2}+h}$,  
$\chi_4'(n-t^2/4)=-\zeta_4^{\epsilon+2h}$ and $J(\eta,\chi_4')=\epsilon+2m\delta \zeta_4$ into \eqref{eq:sizech} and continuing from \eqref{eq:size}, we have $|D_{\ell,h}|=2m^2-m+1$. \qed 
\begin{remark}{\em \label{rem:cond2}
In this remark, we show that there exists a pair $(h,\ell)\in \{0,1,2,3\} \times \{0,1,\ldots,q^2-2\}$ satisfying the condition 
\eqref{eq:cond2}, i.e., 
the set 
\[
\left\{(h,\ell):(q+1) \not | \, \ell, \chi_4'(n)=\zeta_4^{\frac{-\epsilon \delta+1}{2}+h},\chi_4'(n-t^2/4)=-\zeta_4^{\epsilon+2h}\right\}
\]
is nonempty. Let $\Tr_{q^2/q}$ be the trace function from $\F_{q^2}$ to $\F_q$. 
Note that $n-t^2/4=-(\omega^\ell-\omega^{\ell q})^2/4=-\omega^{-(q+1)}\T_{q^2/q}(\omega^{\ell+\frac{q+1}{2}} )^2/4$
is a nonsquare in $\F_q$. Hence,  
\[
\chi_4'(n-t^2/4)=\zeta_4^{-1}\eta(2\omega^{\frac{q-1}{4}}\T_{q^2/q}(\omega^{\ell+\frac{q+1}{2}})).
\] 

Given $\epsilon,\delta\in \{-1,1\}$, set $h$ so that $\frac{-\epsilon \delta+1}{2}+h$ is odd, say, $2d+1$. This is valid whenever $\ell$ is odd since 
$\chi_4'(n)=\chi_4(\omega^\ell)=\zeta_4^{\frac{-\epsilon \delta+1}{2}+h}$. 
Then, the condition $(q+1)\not |\,\ell$ is automatically  satisfied. Furthermore, 
the condition $\chi_4'(n-t^2/4)=-\zeta_4^{\epsilon+2h}$ is equivalent to that  \[
\eta(2\omega^{\frac{q-1}{4}}\Tr_{q^2/q}(\omega^{\ell+\frac{q+1}{2}}))=-\zeta_4^{\epsilon+2h+1}=
-\zeta_4^{\epsilon+2(2d+1-\frac{-\epsilon \delta+1}{2})+1}=
\zeta_4^{\epsilon(1+ \delta)}.
\] 
Therefore, it is enough to see that each of the sets 
\[
T_i=\left\{\omega^\ell\in C_{1}^{(2,q^2)}\,|\,\T_{q^2/q}(\omega^{\ell+\frac{q+1}{2}})\in C_i^{(2,q)}\right\}, \quad i=0,1,
\]
is nonempty. 
The size of each $T_i$ is given by 
\begin{equation}\label{eq:charafun1}
\frac{1}{q}\sum_{a\in \F_q}\sum_{x\in C_{1}^{(2,q^2)}}\sum_{b\in C_i^{(2,q)}}\zeta_p^{\Tr_q(a(\Tr_{q^2/q}(x\omega^{\frac{q+1}{2}})-b))}, 
\end{equation}
where $\Tr_q$ is the trace function from $\F_q$ to the prime field of $\F_q$. 
Let $\psi$ and $\psi'$ be the canonical additive characters of $\F_{q^2}$ and $\F_q$, respectively. Then, 
\begin{align}
\eqref{eq:charafun1}
=&\,\frac{1}{q}\sum_{a\in \F_q}\sum_{x\in C_{1}^{(2,q^2)}}\sum_{b\in C_i^{(2,q)}}\psi(ax\omega^{\frac{q+1}{2}})\psi'(-ab)\nonumber\\
=&\,\frac{1}{q}\sum_{a\in \F_q^\ast}\sum_{x\in C_{1}^{(2,q^2)}}\sum_{b\in C_i^{(2,q)}}\psi(ax\omega^{\frac{q+1}{2}})\psi'(-ab)+\frac{(q-1)(q^2-1)}{4q}. \label{eq:cont1}
\end{align}
Note that $\F_q^\ast \subset C_0^{(2,q^2)}$. Then, by \eqref{eq:Gaussquad} and \eqref{eq:Gaussquad1},  we have 
\begin{align*}
\eqref{eq:cont1}=&\, -\frac{q-1}{2q}\sum_{x\in C_{1}^{(2,q^2)}}\psi(x\omega^{\frac{q+1}{2}})+\frac{(q-1)(q^2-1)}{4q}\\
=&\,-\frac{q-1}{2q}\left(\frac{-1+G_{q^2}(\eta)}{2}\right)+\frac{(q-1)(q^2-1)}{4q}=\frac{q^2-1}{4}. 
\end{align*}
Hence, each $T_i$ is nonempty. 
}\end{remark}

\begin{remark}{\em \label{rem:cond3}
In this remark, we see that the dual of $D_{\ell,h}$ is also a two intersection set with parameters $(2m^2-m;m^2-m,m^2)$ for the block design
obtained from  translations of the set of ``nonsquares'' in $\F_q$ but not for $(P,{\mathcal B})$.
Let $D_{\ell,h}$ be the set defined in Theorem~\ref{thm:twoint}. Let $D_\alpha^\perp=\{s\in \F_q:|D_{\ell,h}\cap (C_0^{(2,q)}+s)|=m^2-m\}$ and    $D_\beta^\perp=\{s\in \F_q:|D_{\ell,h}\cap (C_0^{(2,q)}+s)|=m^2\}$. 
 It is clear that $|D_\alpha^\perp|=2m^2-m$ by \eqref{eq:bsize}. 
By the definitions of $u$ and $\chi_4'$, \eqref{eq:dual} is reformulated as 
\[
N_s=
\begin{cases}
m^2-m& \text{ if  $s\in D_{\ell,h-\epsilon \delta}$,}\\
m^2& \text{ if  $s\in D_{\ell,h+\epsilon \delta}$.}
\end{cases} 
\]
This implies that $D_\alpha^\perp=D_{\ell,h-\epsilon \delta}$ and $D_\beta^\perp=D_{\ell,h+\epsilon \delta}$. 
Then, as in the proof of Theorem~\ref{thm:twoint}, we have  
\begin{align*}
&\, |D_\alpha^\perp\cap  (C_1^{(2,q)}+s)|\\
=
&\, \frac{1}{8}\sum_{x\in \F_q\setminus \{s\}}(-\eta(x-s)+1)\Big(2+\zeta_4^{-h+\epsilon \delta}(1-\zeta_4)\chi_4(1+\omega^\ell x)+\zeta_4^{-3h+3\epsilon \delta}(1+\zeta_4)\chi_4^3(1+\omega^\ell x)\Big)\\
=&\, \frac{1}{8}\Big(
2q-2+\zeta_4^{-h+\epsilon \delta}(1-\zeta_4)(N_{s,1}-N_{s,2})+\overline{\zeta_4^{-h+\epsilon \delta}(1-\zeta_4)(N_{s,1}-N_{s,2})}
\Big) \\
=&\, \begin{cases}
m^2-m& \text{ if  $s\in D_{\ell,h+2}$,}\\
m^2& \text{ if $s\in D_{\ell,h}$.}
\end{cases} 
\end{align*}
This implies that the dual of $D$ forms a two-intersection set for the block 
design 
obtained from translations of $C_1^{(2,q)}$ in $\F_q$.} 
\end{remark}


\section*{Appendix: Upper bounds on excess of conference matrices}
In this appendix, we prove Proposition~\ref{thm:confbound}. 
\begin{proposition}\label{prop:ew}
Let $W$ be a conference matrix of order $n$. 
Then 
$E(W)\leq n\sqrt{n-1}$ holds with equality if and only if 
$W{\bf 1}_n=\sqrt{n-1}{\bf 1}_n$. 
\end{proposition}
\proof 
Write 
\begin{align*}
W{\bf 1}_n=(w_1,w_2\ldots, w_n)^T. 
\end{align*}
Then 
\begin{align*}
\sum_{i=1}^n w_i^2=({\bf 1}_n^T W^T)(W{\bf 1}_n)=(n-1){\bf 1}_n^T {\bf 1}_n=(n-1)n. 
\end{align*}
Thus by Cauchy-Schwartz inequality,
\begin{align*}
E(W)\leq |E(W)|\leq \sum_{i=1}^n |w_i|\leq \sqrt{(\sum_{i=1}^n w_i^2)(\underbrace{1^2+\cdots+1^2}_{n})}=n\sqrt{n-1}.
\end{align*} 
The equality holds if and only if $w_i$ are all equal, that is $w_i=\sqrt{n-1}$ for each $i$. 
Thus, we obtain the assertion. 
\qed
\vspace{0.3cm}

If the equality holds in Proposition~\ref{prop:ew}, then $n-1$ must be a square.  In the next proposition, we improve this upper bound when $n-1$ is a nonsquare. 
\begin{proposition}\label{prop:boundbi}
Let $W$ be a conference matrix of order $n$ with $n-1$ a nonsquare. 
Let $k$ be an odd integer such that $k\leq \sqrt{n-1}<k+2$. 
Then 
\begin{align*}
E(W)\leq \frac{n(k^2+2k+n-1)}{2(k+1)}
\end{align*}
with equality holds if and only if 
$W{\bf 1}_n$ has entries $k,k+2$. 
\end{proposition}
\proof 
Write 
\begin{align*}
W{\bf 1}_n=(w_1,w_2,\ldots, w_n)^T. 
\end{align*}
By the following equation
\begin{align*}
\sum_{i=1}^n(w_i-\sqrt{n-1})^2&=2n(n-1)-2\sqrt{n-1}\sum_{i=1}^n w_i, 
\end{align*}
the value $\sum_{i=1}^n w_i$ takes maximum if and only if the value $\sum_{i=1}^n(w_i-\sqrt{n-1})^2$ takes minimum. 
The latter occur only if each $w_i$ is either $k$ or $k+2$. Here, noting that 
the sum of entries in each row has the same parity with $n-1$, $w_i=k+1$ is impossible. 
In this case, 
\begin{align*}
W{\bf 1}_n=\begin{pmatrix} k{\bf 1}_a\\ (k+2) {\bf 1}_{n-a} \end{pmatrix},  
\end{align*}
and the number $a$ of $k$ in $W{\bf 1}_n$ is determined as $a=\frac{n((k+2)^2-(n-1))}{4(k+1)}$ by $\sum_{i-1}^n w_i^2=n(n-1)$. 
Then 
\begin{align*}
E(W)&=\sum_{i=1}^n w_i=\frac{1}{2\sqrt{n-1}}(2n(n-1)-\sum_{i-1}^n(w_i-\sqrt{n-1})^2)\\
&\leq n\sqrt{n-1}-\frac{1}{2\sqrt{n-1}}(a (k-\sqrt{n-1})^2+(n-a) (k+2-\sqrt{n-1})^2)
=\frac{n(k^2+2k+n-1)}{2(k+1)}
\end{align*}
with  equality holds if and only if  $W{\bf 1}_n$ has only two entries $k,k+2$. 
\qed
\vspace{0.3cm}

Combining Propositions~\ref{prop:ew} and \ref{prop:boundbi}, we have Proposition~\ref{thm:confbound}. Note that  Propositions~\ref{prop:ew} and \ref{prop:boundbi} are generalizable for 
excess of general {\it weighing matrices}. 
\end{document}